\documentclass[11pt]{amsart}
\usepackage{amsmath,amssymb} 
\usepackage{graphicx}
\usepackage[utf8]{inputenc} 
\usepackage[font=small,labelfont=bf]{caption}

\begin{document}

\newtheorem{thm}{Theorem}[section]
\newtheorem{lem}[thm]{Lemma}
\newtheorem{prop}[thm]{Proposition}
\newtheorem{coro}[thm]{Corollary}
\newtheorem{defn}[thm]{Definition}
\newtheorem*{remark}{Remark}

\numberwithin{equation}{section}

\newcommand{\Z}{{\mathbb Z}} %cph changed from \mathbf
\newcommand{\Q}{{\mathbb Q}}
\newcommand{\PP}{{\mathbb P}}
\newcommand{\R}{{\mathbb R}}
\newcommand{\C}{{\mathbb C}}
\newcommand{\N}{{\mathbb N}}
\newcommand{\FF}{{\mathbb F}}
\newcommand{\T}{{\mathbb T}}
\newcommand{\fq}{\mathbb{F}_q}
\newcommand{\IS}{{\mathbb S}}

\newcommand{\fixmehidden}[1]{}

\def\scrA{{\mathcal A}}
\def\cB{{\mathcal B}}
\def\Eps{{\mathcal E}}
\def\cI{{\mathcal I}}
\def\scrD{{\mathcal D}}
\def\cF{{\mathcal F}}
\def\cL{{\mathcal L}}
\def\cM{{\mathcal M}}
\def\cN{{\mathcal N}}
\def\cP{{\mathcal P}}
\def\scrC{{\mathcal C}}
\def\scrR{{\mathcal R}}
\def\scrS{{\mathcal S}}

\newcommand{\rmk}[1]{\footnote{{\bf Comment:} #1}}

\renewcommand{\mod}{\;\operatorname{mod}}
\newcommand{\ord}{\operatorname{ord}}
\newcommand{\TT}{\mathbb{T}}
\renewcommand{\i}{{\mathrm{i}}}
\renewcommand{\d}{{\mathrm{d}}}
\renewcommand{\^}{\widehat}
\newcommand{\HH}{\mathbb H}
\newcommand{\Vol}{\operatorname{vol}}
\newcommand{\area}{\operatorname{area}}
\newcommand{\tr}{\operatorname{tr}}
\newcommand{\norm}{\mathcal N} % norm =(\frac{ n+\sqrt{n^2-4}} 2)^2
\newcommand{\intinf}{\int_{-\infty}^\infty}
\newcommand{\ave}[1]{\left\langle#1\right\rangle} %  average
\newcommand{\E}{\mathbb E}
\newcommand{\Var}{\operatorname{Var}}
\newcommand{\Cov}{\operatorname{Cov}}
\newcommand{\Prob}{\operatorname{Prob}}
\newcommand{\sym}{\operatorname{Sym}}
\newcommand{\disc}{\operatorname{disc}}
\newcommand{\CA}{{\mathcal C}_A}
\newcommand{\cond}{\operatorname{cond}} % conductor
\newcommand{\lcm}{\operatorname{lcm}}
\newcommand{\Kl}{\operatorname{Kl}} %Kloosterman sum
\newcommand{\leg}[2]{\left( \frac{#1}{#2} \right)}  % Legendre symbol
\newcommand{\id}{\operatorname{id}}
\newcommand{\beq}{\begin{equation}}
\newcommand{\eeq}{\end{equation}}
\newcommand{\bsp}{\begin{split}}
\newcommand{\esp}{\end{split}}
\newcommand{\bra}{\left\langle}
\newcommand{\ket}{\right\rangle}
\newcommand{\diam}{\operatorname{diam}}
\newcommand{\supp}{\operatorname{supp}}
\newcommand{\dist}{\operatorname{dist}}
\newcommand{\sgn}{\operatorname{sgn}}
\newcommand{\inte}{\operatorname{int}}
\newcommand{\ind}{\operatorname{ind}}
\newcommand{\Spec}{\operatorname{Spec}}
\newcommand{\sumstar}{\sideset \and^{*} \to \sum}

\newcommand{\LL}{\mathcal L} %L-function of u
\newcommand{\sumf}{\sum^\flat}
\newcommand{\Hgev}{\mathcal H_{2g+2,q}}
\newcommand{\USp}{\operatorname{USp}}
\newcommand{\conv}{*}
\newcommand{\CF}{c_0} % Fejer constant
\newcommand{\kerp}{\mathcal K}

\newcommand{\gp}{\operatorname{gp}}
\newcommand{\Area}{\operatorname{Area}}

\title[]{Unique continuation at infinity for potentials with arbitrary radial growth}
%Landis-type theorems for potentials with arbitrary radial growth on $\R^n$} 
  
\author{Henrik Uebersch\"ar}
\address{Sorbonne Universit\'e and Universit\'e Paris Cit\'e, CNRS, IMJ-PRG, F-75005 Paris, France.}
\email{henrik.ueberschar@imj-prg.fr}
\date{\today}

\date{\today}
\maketitle

\begin{abstract}
Let $G$ be any given continuous positive function on $\R_+$. Let $V$ be radial with $|V(x)|\leq G(|x|)$. We prove a Landis-type theorem for any real-valued solution of $\Delta u=Vu$ on $\R^n$. We construct a decay threshold $e^{-g(r)}$, where $g$ is a strictly increasing function which can be computed explicitly in terms of $G$. 
Under suitable assumptions the exponent in the decay threshold is proportional to the Agmon distance associated with $G$.  
%For the special case of potentials with polynomial growth our decay threshold matches the sharp threshold recently proved by Davey on the plane. 
%Moreover, we show that there exist functions $G$ of arbitrarily rapid growth which still admit an exponential decay threshold provided that $G$ is allowed to vary sufficiently strongly. 
%Finally, our result implies a quantitative version of the strong Landis conjecture for the case of radial potentials.
\end{abstract}

\section{Introduction}
\subsection{Background}
Let $V$ be a bounded, measurable function on $\R^n$.
A well-known conjecture of Landis \cite{Lan53,Lan71} asserts that there exists a constant $c>0$ large enough such that any real-valued solution of $\Delta u=Vu$ on $\R^n$ which satisfies the bound $|u(x)|\leq e^{-c|x|}$ must be trivial. Meshkov showed \cite{Me92} that an analogous result does not hold for complex-valued solutions, where the decay threshold is of order $e^{-c|x|^{4/3}}$. The real-valued case of the conjecture distinguishes between its {\em qualitative} version (i.e. any solution satisfying the decay on $\R^n$ is trivial) as well as a {\em quantitative} version of the conjecture (explicit lower bounds for a nontrivial solution). The original conjecture with exponential decay is sometimes referred to as the {\em strong} Landis conjecture, as opposed to weaker versions which allow for a small loss in the exponent.

A weak version of Landis' conjecture was proved by Logunov, Malinnikova, Nadirashvilli and Nazarov \cite{LMNN24} in dimension $n=2$. %We note that the authors in \cite{LMNN24} proved both a qualitative and a quantitative version of the conjecture. 
In dimension $n\geq 3$, the conjecture was recently shown to be false by Frank and Ivanisvili \cite{Fr26}. Their counterexamples match Meshkov's decay treshold. However, in the case of radial potentials the qualitative version of the conjecture holds in any dimension, as was shown by Rossi in 2018 \cite{Ros21}.

A natural generalization of Landis' conjecture is to extend the problem to potentials which may grow at infinity, and one expects that the decay threshold ought to depend on the growth of the potential. For example, a sharp decay threshold has been established by Davey for the case of potentials with polynomial growth \cite{Dav24} utilising the methods of \cite{LMNN24}.

\subsection{Statement of results}
Let the potential $V=V(|x|)$ be radial. We consider a real-valued (not necessarily radial) solution $u$ of $\Delta u=Vu$ on $\R^n$. We allow $|V(r)|$ to grow arbitrarily rapidly, as $r\to+\infty$.

We have two main results: the first applies to any continuous $G\geq0$ without assumptions on its variation (see Thm. \ref{main}); the second result assumes $G\in C^2$, $G\geq c$ for some constant $c>0$, and imposes conditions on the variation of $G$ (see Thm. \ref{main2}). 
%Our first result proves a decay threshold for a solution $u$, where the potential satisfies a growth condition $|V|\leq G$, where $G\in C^0(\R_+,\R_+)$ is any given growth function. 
%We prove the following quantitative Landis-type theorem. The proof makes crucial use of the real-valuedness of the solution. 
\begin{thm}\label{main}
Let $G\in C^0(\R_+,\R_+)$ be nonzero and $G\geq0$. Define
\beq
\beta(r):=(n+1)r+\int_0^r G.
\eeq
Let $u$ be a real-valued solution of $\Delta u=Vu$ on $\R^n$, where $V$ is radial and satisfies $|V|\leq G$. 
If $u(0)\neq 0$, then there exists a constant $c>0$ and a sequence of points $(x_i)_{i=0}^{+\infty}$ in $\R^n$ such that $\lim_{i\to+\infty}|x_i|=+\infty$ and $|u(x_i)|\geq ce^{-\beta(|x_i|)}$.
% i.e. negation of $\lim_{|x|\to+\infty}|u(x)|e^{-\beta(|x|)}=0$
\end{thm}

Bourgain-Kenig \cite{BoKe05} raised the following question about the decay properties of a real-valued solution $u$ when $|V|\leq 1$: which estimate can be given on
\begin{equation}
B_u:=\liminf_{|x|\to+\infty} \frac{\log\log |u(x)|^{-1}}{\log |x|}?
\end{equation}

Logunov, Malinnikova, Nadirashvilli and Nazarov gave a quantitative lower bound of order $\exp(-c r(\log r)^{3/2})$ on the supremum of $|u|$ on a unit ball centered at $|x|=r/2>2$ for any solution satisfying $|u(0)|=1=\sup_{B(0,2r)}|u|$. This is sufficient to show that $B_u\leq 1$. 

However, due to the log-loss in the exponent, their result is not sufficient to bound the quantity
\begin{equation}
A_u:=\liminf_{|x|\to+\infty} \frac{\log |u(x)|^{-1}}{|x|}.
\end{equation}

As a consequence of Theorem \ref{main} we can prove the following estimate.
\begin{coro} \label{cor1} 
Let $V\in L^\infty(\R^n)$ be radial and satisfy $|V|\leq 1$.
Let $u$ be a real-valued solution of $\Delta u=Vu$ on $\R^n$ which satisfies $u(0)\neq0$. Then $A_u\leq n+2$.
\end{coro}

Another application of Theorem \ref{main} is the construction of potentials of arbitrary radial growth such that solutions still satisfy an exponential lower bound along a sequence of points.
\begin{coro}\label{cor2}
Let $f:\R_+\to\R_+$ be any strictly increasing function such that $(f(n)^{-1})_{n\geq1}$ is summable. Then there exists a radial potential such that $|V(x)|=f(n)$ on every circle $|x|=n$, for $n\geq 1$, and any real-valued solution of $\Delta u=Vu$, with $u(0)\neq0$, admits a quantitative exponential lower bound in the sense that there exists a constant $K_f>0$ and  a sequence of points $(x_i)_{i=0}^{+\infty}$ such that $\lim_{i\to+\infty}|x_i|=+\infty$ and $|u(x_i)|\geq K_f e^{-(n+2)|x_i|}$.
\end{coro}

Under additional assumptions on the function $G$ we are able to prove the following result.  We note that the exponent function in the theorem below seems to be closely related to the Agmon metric \cite{Ag79} which appears as the exponent function in exponential decay estimates of eigenfunctions for potentials $V$ s. t. $V(x)\to+\infty$, as $|x|\to+\infty$
\begin{thm}\label{main2}
Let $G\in C^2(\R_+,\R_+)$, and assume there exists $c>0$ s. t. $G\geq c$. We define 
\beq\label{exp}
g(r)=\int_0^r \sqrt{G(s)}ds.
\eeq
Assume that $f=g^{-1}:\R_+\mapsto\R_+$ satisfies the condition that $|f''/f'|$ and $|f'''/f'|$ are bounded on $[1,+\infty)$.

Let $u$ be a real-valued solution of $\Delta u=Vu$ on $\R^n$, where $V$ is radial and satisfies $|V|\leq G$. 
If $u(0)\neq 0$, then there exists $c>0$ and a sequence $(x_j)$ in $\R^n$, $|x_j|\to+\infty$, s. t. $|u(x_j)|\geq c e^{-\beta g(|x_j|)}$, where $\beta>0$ is an absolute constant which depends only on the function $G$.
\end{thm}
We note that the theorem above yields exactly the type of decay which appears in the sharp cutoff derived by Davey in \cite{Dav24} for the case of potentials of polynomial growth $G(r)=(1+r)^N$, $N>0$. Namely $g(r)\sim\tfrac{2}{N+2}r^{N/2+1}$.\\

{\bf Conflict of interest statement.}

The author states that there is no conflict of interest.\\

{\bf Data availability statement.}

No data was produced as part of this research.

\section{Proofs}

\subsection{Proof of Theorem \ref{main}}
We will prove the claim by contradiction. So let us assume that $$\lim_{r\to+\infty}e^{\beta(r)}\sup_{|x|=r}|u(x)|=0.$$
We denote by $u(r,\theta)$ the solution written in spherical coordinates $(r,\theta)\in \R_+\times \IS^{n-1}$. Let us define the function
\beq
w(r):=\frac{1}{\Vol(\IS^{n-1})}\int_{\IS^{n-1}} u(r,\theta).
\eeq

Since $u$ satisfies the Schr\"odinger equation in spherical coordinates
\beq
\partial_r^2 u+(n-1)r^{-1}\partial_r u+r^{-2}\Delta_\theta u =V(r)u,
\eeq
where $\Delta_\theta$ denotes the Laplace-Beltrami operator on $\IS^{n-1}$,
we have
\beq
w''(r)=\frac{1}{\Vol(\IS^{n-1})}\int_{\IS^{n-1}} \partial_r^2 u(r,\theta)
=(1-n)r^{-1}w'(r)+V(r)w(r),
\eeq
where we used $\int_{\IS^{n-1}}\Delta_\theta u =0$.

Our assumption implies that $w$ satisfies $|w(r)| e^{\beta(r)}\to 0$, as $r\to+\infty$. Note that $w$ does not vanish identically on $[1,+\infty)$: To see this, recall our assumption that $u(0)\neq0$ which implies, by continuity of $w$, that $w\neq 0$ on an interval $(0,\epsilon]$ for some $\epsilon>0$. 

Now, $w$ solves 
\begin{equation}\label{ode}
w''(r)=(1-n)r^{-1}w'(r)+V(r)w(r)
\end{equation} 
on $[\epsilon/2,+\infty)$, a second order ODE with bounded coefficients. 
%Therefore, $w$ satisfies Carleman's unique continuation theorem \cite{Car39,Ler19}.  %This is much easier...
Assume, for a contradiction, that $w$ vanishes identically on $[1,+\infty)$. Then, by unique continuation, $w$ vanishes identically on $[\epsilon/2,+\infty)$ which contradicts the nonvanishing of $w$ on $[0,\epsilon]$.

Since $w$ does not vanish identically on $[1,+\infty)$, we have that $h(r):=|w(r)|e^{\beta(r)}$ does not vanish identically on $[1,+\infty)$. Moreover, we note that $h$ is continuous and $h(r)\to 0$, as $r\to+\infty$. This means $C_0:=\sup_{[1,+\infty)} h>0$ is attained on $[1,+\infty)$. Hence, there exists $r_0\geq 1$ s. t. $h(r_0)=C_0$ and, therefore, $|w(r_0)|=C_0 e^{-\beta(r_0)}$ and $|w(r)|\leq C_0 e^{-\beta(r)}$ for any $r\geq 1$.

Moreover, we obtain a lower bound on $|w'(r_0)|$: without loss of generality, assume $w(r_0)>0$. Due to continuity of $w$ there exists $\delta>0$ s. t. $w(r)>0$ for any $r_0<r<r_0+\delta$. We then have
$$\frac{w(r)-w(r_0)}{r-r_0}=\frac{w(r)-C_0 e^{-\beta(r_0)}}{r-r_0}\leq C_0\frac{e^{-\beta(r)}-e^{-\beta(r_0)}}{r-r_0},$$
where we used $w(r)=|w(r)|\leq C_0e^{-\beta(r)}$ and $r>r_0$.

If we pass to the limit $r\to r_0^+$, we obtain $w'(r_0)\leq -C_0\beta'(r_0)e^{-\beta(r_0)}$, and, thus, $|w'(r_0)|\geq C_0\beta'(r_0)e^{-\beta(r_0)}$.

In summary, there exists $C_0>0$ and $r_0\geq 1$ s. t. $|w(r_0)|=C_0e^{-\beta(r_0)}$, $|w'(r_0)|\geq C_0\beta'(r_0)e^{-\beta(r_0)}$, and $|w(r)|\leq C_0e^{-\beta(r)}$ for any $r\geq1$.

For $T>2r_0$ we define a smooth cutoff $\chi:[r_0,+\infty)\to[0,1]$
\beq
\chi(r)=
\begin{cases}
1, \quad r_0\leq r < T/2,\\
\\
0, \quad r>T.
\end{cases}
\eeq
and $\chi$ is chosen on $[T/2,T]$ s. t. $\chi(T/2)=1$, $\chi(T)=0$ and $\chi$ is smooth and decreasing with $|\chi'|\leq cT^{-1}$, $|\chi''|\leq c'T^{-2}$.

Let us integrate \eqref{ode} on the interval $[r_0,T]$. \\
%(where we assume $w'(r)\to 0$, as $r\to+\infty$)
\beq
\begin{split}
w'(r_0)=&-\int_{r_0}^{T}w''(r)\chi(r)-\int_{r_0}^T w'(r)\chi'(r)\\
=& (n-1)\int_{r_0}^{T}r^{-1}w'(r)\chi(r)-\int_{r_0}^{T}V(r)w(r)\chi(r)-\int_{r_0}^T w'(r)\chi'(r).
\end{split}
\eeq

We integrate by parts to see
\beq
\begin{split}
w'(r_0)=& -(n-1)r_0^{-1}w(r_0)\\
&+\int_{r_0}^{T}w(r)((n-1)(r^{-2}\chi(r)-r^{-1}\chi'(r))+\chi''(r))\\
&-\int_{r_0}^{T}V(r)w(r)\chi(r)
\end{split}
\eeq

From this we obtain the inequality
\beq
\begin{split}
|w'(r_0)|\leq& (n-1)r_0^{-1}|w(r_0)|\\
&+(n-1)\int_{r_0}^{+\infty}(r^{-2}+cr^{-1}T^{-1}+c'T^{-2})|w(r)|\\
&+\int_{r_0}^{+\infty}|V(r)||w(r)|
\end{split}
\eeq
which is valid for any $T$ large enough. 

If we send $T\to +\infty$, we obtain the inequality
\beq
\begin{split}
|w'(r_0)|\leq (n-1)r_0^{-1}|w(r_0)|+&(n-1)\int_{r_0}^{+\infty}r^{-2}|w(r)|\\
&+\int_{r_0}^{+\infty}|V(r)||w(r)|.
\end{split}
\eeq

%{\bf
%NOTE: Using the point-wise gradient estimate from the $L^\infty$ paper we can always get a bound of the type $|w'(r)|\lesssim_V e^{-\beta(r-1)}$, where the implied constant depends on local $L^\infty$ bounds on a number of derivatives of $V$. As long as the potential doesn't oscillate very fast, then this should work. But it requires $V\in C^k$, where $k=k(n)$, which goes to infinity, as $n\to+\infty$. 

%In fact, the bound we would use is of the form $$|\partial_r u(r,\theta)|\lesssim_V \|u\|_{L^2(B_1(r,\theta))},$$ where $B_1(r,\theta)$ denotes a unit ball centered on the point $(r,\theta)$ and the implied constant depends on local $L^\infty$ bounds on a number of derivatives of $V$.
%}

If we use the estimates above, we obtain
\beq
\begin{split}
C_0\beta'(r_0)e^{-\beta(r_0)}\leq |w'(r_0)|\leq &(n-1)C_0e^{-\beta(r_0)}\\
&+C_0\int_{r_0}^{+\infty}(n-1+|V(r)|)e^{-\beta(r)}.
\end{split}
\eeq

Thus, we have
\beq\label{ineq}
\beta'(r_0)\leq n-1 + e^{\beta(r_0)}\int_{r_0}^{+\infty}(n-1+|V(r)|)e^{-\beta(r)}.
\eeq

We now recall that, by definition of the function $\beta$, we have for any $r\geq 1$:
\beq
\beta'(r)= n+1+G(r).
\eeq
%n-1+|V(r)|\leq 
Moreover, we estimate (recall $|V(r)|\leq G(r)$)
\begin{equation}
\begin{split}
\int_{r_0}^{+\infty}(n-1+|V(r)|)e^{-\beta(r)}
&\leq \int_{r_0}^{+\infty}(n+1+G(r))e^{-\beta(r)}\\
&=\int_{r_0}^{+\infty}\beta'(r)e^{-\beta(r)}=e^{-\beta(r_0)}
\end{split}
\end{equation}
and, by substituting this inequality in \eqref{ineq}, we obtain
$$\beta'(r_0)\leq n-1 + e^{\beta(r_0)}\int_{r_0}^{+\infty}\beta'(r)e^{-\beta(r)}=n.$$

And, because $G\geq0$, we arrive at a contradiction:
$$n+1\leq n+1+G(r_0)=\beta'(r_0)\leq n$$

Hence, we have shown that there exists $\epsilon>0$ and a sequence of radii $r_j\nearrow +\infty$ s. t. $e^{\beta(r_j)}\sup_{|x|=r_j}|u(x)|\geq \epsilon.$ On each sphere of radius $r_j$ the supremum of $|u|$ is attained. So there exists $x_j\in\R^n$, $|x_j|=r_j$, s. t. $|u(x_j)|\geq\epsilon e^{-\beta(|x_j|)}$.

\subsection{Proof of Corollary \ref{cor1}} In the case $G=1$, Theorem \ref{main} gives $\beta(r)=(n+2)r$. Hence, any real-valued solution $u$ of $\Delta u=Vu$, $|V|\leq 1$, which satisfies $u(0)\neq 0$, admits a lower bound $|u(x_i)|\geq ce^{-(n+2)|x_i|}$ along a sequence $(x_i)_{i=0}^{+\infty}$, $|x_i|\to+\infty$. Hence, for $|x_i|\geq c^{-1}$ we have $|u(x_i)|\geq ce^{-(n+2)|x_i|}\geq |x_i|^{-1}e^{-(n+2)|x_i|}$. This yields $|u(x_i)|^{-1}\leq |x_i|e^{(n+2)|x_i|}$ and, thus, $|x_i|^{-1}\log |u(x_i)|^{-1} \leq n+2+|x_i|^{-1}\log|x_i|$ which implies the result.\\

\subsection{Proof of Corllary \ref{cor2}}
It is worth noting that the growth of the function $G$ is, in fact, irrelevant for the decay threshold $e^{-\beta(r)}$, because we only require $G$ to be continuous -- not monotone. One can, therefore, easily construct a continuous function $G$ of arbitrary growth along a subsequence $(x_j)_{j=0}^{+\infty}\subset\R_+$, s. t. $\int_0^{+\infty} G<+\infty$, so that the decay threshold is still exponential. For example let $(r_n)_{n\geq 1}$ be a decreasing sequence such that for all $n\geq 1$: $r_n>0$ and $r_n\searrow 0$. Moreover, assume that $(r_n)$ satisfies
$$\sum_{n=1}^{+\infty} r_n<+\infty.$$

Fix $\chi\in C^\infty_c(\R)$ with $0\leq\chi\leq 1$ and $\supp\chi\subset[-\tfrac{1}{2},\tfrac{1}{2}]$ s. t. $\chi=1$ on $[-\tfrac{1}{4},\tfrac{1}{4}]$. We construct
\beq
G(s)=\sum_{n=1}^{+\infty} r_n^{-1}\chi\left(\frac{s-n}{r_n^2}\right)
\eeq
which satisfies $$\int_0^{+\infty} G= \hat{\chi}(0)\sum_{n=1}^{+\infty}r_n, \quad\text{where}\quad \hat{\chi}(0)=\int_{\R}\chi.$$

Let $f:\R_+\mapsto\R_+$ be any function which is sufficiently rapidly increasing to ensure that $f(n)^{-1}$ is summable. Define $r_n=f(n)^{-1}$.

Then we have $G(n)=f(n)$ and $$C_f:=\int_0^{+\infty}G=\hat{\chi}(0)\sum_{n\geq 1} f(n)^{-1}.$$

We apply Theorem \ref{main} and have the following estimate for the function $\beta$:
$$\beta(r)=(n+1)r+\int_0^r G\leq (n+1)r+C_f$$

Thus, if the potential satisfies the growth condition $|V(x)|\leq G(|x|)$ for any $x\in\R^n$, then for any real-valued solution $u$ of $\Delta u=Vu$ which satisfies $u(0)\neq0$ we have the following: there exists a constant $c>0$ and a sequence of points $(x_i)$ s. t.  $|x_i|\to+\infty$ and for $e^{|x_i|}\geq c^{-1}$ we have $$|u(x_i)|\geq c e^{-C_f} e^{-(n+1)|x_i|}\geq e^{-C_f}e^{-(n+2)|x_i|}.$$

\subsection{Proof of Theorem \ref{main2}}

Define a strictly increasing function $g:\R_+\mapsto\R_+$ by $$g(r)=\int_0^r\sqrt{G(s)}ds.$$
We will prove the claim by contradiction. So let us assume that $$\lim_{r\to+\infty}e^{\beta g(r)}\sup_{|x|=r}|u(x)|=0.$$
We denote by $u(r,\theta)$ the solution written in spherical coordinates $(r,\theta)\in \R_+\times \IS^{n-1}$. Let us define the function
\beq
w(r):=\frac{1}{\Vol(\IS^{n-1})}\int_{\IS^{n-1}} u(r,\theta)
\eeq
and denote $v(r)=w(f(r))$, where $f=g^{-1}$.

Our assumption implies that $v$ satisfies $|v(r)| e^{\beta r}\to 0$, as $r\to+\infty$. Note that $v$ does not vanish identically on $[1,+\infty)$: To see this, recall our assumption that $u(0)\neq0$ which implies, by continuity of $w$, that $w\neq 0$ on an interval $(0,\epsilon]$ for some $\epsilon>0$. 

Now, $w$ solves 
\begin{equation}\label{ode2}
w''(r)=(1-n)r^{-1}w'(r)+V(r)w(r).
\end{equation} 

We compute $v'(r)=f'(r)w'(f(r))$ and $$v''(r)=f''(r)w'(f(r))+f'(r)^2 w''(f(r))$$
where we substitute \eqref{ode2} and $w'(f(r))=(f'(r))^{-1}v'(r)$ to obtain
$$v''(r)=\frac{f''(r)}{f'(r)}v'(r)+f'(r)^2[ (1-n)f(r)^{-1}w'(f(r))+V(f(r))w(f(r))]$$
and, by substituting $w'(f(r))=(f'(r))^{-1}v'(r)$ again, we obtain
\begin{equation}\label{ode3}
v''(r)=\left(\frac{f''(r)}{f'(r)}+(1-n)\frac{f'(r)}{f(r)}\right)v'(r)+f'(r)^2V(f(r))v(r)
\end{equation}

In particular, $v$ solves the above second order ODE with bounded coefficients on $[\epsilon/2,+\infty)$.
Assume, for a contradiction, that $v$ vanishes identically on $[1,+\infty)$. Then, by unique continuation, $v$ vanishes identically on $[\epsilon/2,+\infty)$ which contradicts the nonvanishing of $v$ on $[0,\epsilon]$.

Since $v$ does not vanish identically on $[1,+\infty)$, we have that $h(r):=|v(r)|e^{\beta r}$ does not vanish identically on $[1,+\infty)$. Moreover, we note that $h$ is continuous and $h(r)\to 0$, as $r\to+\infty$. This means $C_0:=\sup_{[1,+\infty)} h>0$ is attained on $[1,+\infty)$. Hence, there exists $r_0\geq 1$ s. t. $h(r_0)=C_0$ and, therefore, $|v(r_0)|=C_0 e^{-\beta r_0}$ and $|v(r)|\leq C_0 e^{-\beta r}$ for any $r\geq 1$.

Moreover, we obtain a lower bound on $|v'(r_0)|$: without loss of generality, assume $v(r_0)>0$. Due to continuity of $v$ there exists $\delta>0$ s. t. $v(r)>0$ for any $r_0<r<r_0+\delta$. We then have
$$\frac{v(r)-v(r_0)}{r-r_0}=\frac{v(r)-C_0 e^{-\beta(r_0)}}{r-r_0}\leq C_0\frac{e^{-\beta r}-e^{-\beta r_0}}{r-r_0},$$
where we used $v(r)=|v(r)|\leq C_0e^{-\beta r}$ and $r>r_0$.

If we pass to the limit $r\to r_0^+$, we obtain $v'(r_0)\leq -C_0\beta e^{-\beta r_0}$, and, thus, $|v'(r_0)|\geq C_0\beta e^{-\beta r_0}$.

In summary, there exists $C_0>0$ and $r_0\geq 1$ s. t. $|v(r_0)|=C_0e^{-\beta r_0}$, $|v'(r_0)|\geq C_0\beta e^{-\beta r_0}$, and $|v(r)|\leq C_0e^{-\beta r}$ for any $r\geq1$.

For $T>2r_0$ we define a smooth cutoff $\chi:[r_0,+\infty)\to[0,1]$
\beq
\chi(r)=
\begin{cases}
1, \quad r_0\leq r < T/2,\\
\\
0, \quad r>T.
\end{cases}
\eeq
and $\chi$ is chosen on $[T/2,T]$ s. t. $\chi(T/2)=1$, $\chi(T)=0$ and $\chi$ is smooth and decreasing with $|\chi'|\leq c_3T^{-1}$, $|\chi''|\leq c_4T^{-2}$.

Let us integrate \eqref{ode3} on the interval $[r_0,T]$
\beq
\begin{split}
v'(r_0)=&-\int_{r_0}^{T}v''(r)\chi(r)-\int_{r_0}^T v'(r)\chi'(r)\\
=&-\int_{r_0}^{T}\alpha(r)v'(r)\chi(r)-\int_{r_0}^{T}f'(r)^2V(f(r))v(r)\chi(r)-\int_{r_0}^T v'(r)\chi'(r),
\end{split}
\eeq
where we denote
$$\alpha(r)=\frac{f''(r)}{f'(r)}+(1-n)\frac{f'(r)}{f(r)}.$$

We have 
$$\alpha'(r)=\frac{f'''(r)f'(r)-f''(r)^2}{f'(r)^2}+(1-n)\frac{f''(r)f(r)-f'(r)^2}{f(r)^2}.$$

We integrate by parts to see
\beq
\begin{split}
v'(r_0)=& \alpha(r_0)v(r_0)\\
&+\int_{r_0}^{T}(\alpha'(r)\chi(r)+\alpha(r)\chi'(r)+\chi''(r))v(r)\\
&-\int_{r_0}^{T}f'(r)^2V(f(r))v(r)\chi(r)
\end{split}
\eeq

Recall our assumption that $|f''(r)|/|f'(r)|$, $|f'''(r)|/|f'(r)|$ are bounded for any $r\geq 1$. Moreover, our assumption $G\geq c>0$ implies
$$|f'(r)|=|(G(f(r)))^{-1/2}|\leq c^{-1/2}.$$

This implies that there exist constants $c_1,c_2>0$ s. t. for any $r\geq1$
\beq
|\alpha(r)|\leq c_1, \quad |\alpha'(r)|\leq c_2.
\eeq
From this we obtain the inequality
\beq
\begin{split}
|v'(r_0)|\leq& c_1|v(r_0)|\\
&+\int_{r_0}^{+\infty}(c_2+c_1c_3T^{-1}+c_4T^{-2})|v(r)|\\
&+\int_{r_0}^{+\infty}f'(r)^2|V(f(r))||v(r)|
\end{split}
\eeq
which is valid for any $T$ large enough. 

If we send $T\to +\infty$, we obtain the inequality
\beq
|v'(r_0)|\leq c_1|v(r_0)|+c_2\int_{r_0}^{+\infty}|v(r)|+\int_{r_0}^{+\infty}f'(r)^2|V(f(r))||v(r)|.
\eeq

We note that we have the estimate
\beq
f'(r)^2|V(f(r))|\leq f'(r)^2G(f(r))=1
\eeq
where we recall that $f^{-1}=g$ and $g'(s)=\sqrt{G(s)}$.

If we use the estimates above, we obtain
\beq
C_0\beta e^{-\beta r_0}\leq |v'(r_0)|\leq c_1C_0e^{-\beta r_0}+(1+c_2)C_0\int_{r_0}^{+\infty}e^{-\beta r}.
\eeq

If we evaluate the integral, we obtain
$$C_0\beta e^{-\beta r_0}\leq |v'(r_0)|\leq c_1C_0e^{-\beta r_0}+(1+c_2)C_0\beta^{-1}e^{-\beta r_0}$$
and upon multiplication by $C_0^{-1}e^{\beta r_0}$ we have
$$\beta \leq c_1+(1+c_2)\beta^{-1}$$
which is a contradiction for $\beta$ large enough.

Hence, we have shown that there exists $\epsilon>0$ and a sequence of radii $r_j\nearrow +\infty$ s. t. $e^{\beta g(r_j)}\sup_{|x|=r_j}|u(x)|\geq \epsilon.$ On each sphere of radius $r_j$ the supremum of $|u|$ is attained. So there exists $x_j\in\R^n$, $|x_j|=r_j$, s. t. $|u(x_j)|\geq\epsilon e^{-\beta g(|x_j|)}$.

\end{document}